\documentstyle[12pt]{article}
\newcommand{\card}{\mathop{\rm card}\limits}
\newcommand{\const}{\mathop{\rm const}\limits}

\newcommand{\mod}{\mathop{\rm mod}\limits}
\newcommand{\Var}{\mathop{\rm Var}\limits}

\newcommand{\supp}{\mathop{\rm supp}\limits}

\begin{document}
\begin{center}

{\bf MOMENT AND TAIL INEQUALITIES } \\

\vspace{3mm}

{\bf FOR POLYNOMIAL MARTINGALES.}\\

\vspace{4mm}

{\bf The case of heavy tails.}\\

\vspace{3mm}
{\bf E. Ostrovsky and L. Sirota}\\

\vspace{2mm}

{\it Department of Mathematics and Statistics, Bar-Ilan University,
59200, Ramat Gan.}\\
e \ - \ mails: eugostrovsky@list.ru; \ sirota3@bezeqint.net \\

\vspace{4mm}

                 {\it Abstract.}\\
\end{center}

 In this paper we obtain the non-asymptotic exact moment and tails
estimates for polynomial on martingale differences. \par
 We give also some examples on order to show the exactness of obtained
results. \par
 \vspace{3mm}

 \par

 \vspace{3mm}

\par

\vspace{3mm}

\section{Introduction. Notations. Statement of problem.}

\vspace{3mm}

 Let $ (\Omega,F,{\bf P} ) $ be a probabilistic space, $ \Omega = \{\omega\},$
$  \xi(i,1), \xi(i,2), \ldots, \xi(i,d) $  be a family of centered:
$ {\bf E} \xi(i,m) = 0, \ m=1,2,\ldots,d $ martingale differences relative some
filtration $ \{ F(i) \}: $

$$
F(0) = \{ \emptyset, \Omega \}, F(i) \subset F(i+1) \subset F:
$$

 $  \forall k=0,1,\ldots, i-1 \ \Rightarrow $
$$
{\bf E}\xi(i,m)/F(k)=0; \ {\bf E} \xi(i,m)/F(i)= \xi(i,m) \ (\mod {\bf P}), \eqno(1.1)
$$
"martingale version".\par

\vspace{3mm}

 In particular, the random {\it vectors} $ \vec{\xi}(m) = \{ \xi(i,m) \}, i=1,2,\ldots,n $
for all the values $ n=1,2,\ldots $ may be independent; we will call this possibility
" vector independent version". \par

 A second case: the random vectors $ \vec{\xi}(m) = \{ \xi(i,m) \}, i=1,2,\ldots,n $
 may be dependent, but {\it inside} of this vectors the random variables $ \{ \xi(i,m) \}, \
 i=1,2,\ldots $  are independent for all the values $ m; \ m=1,2,\ldots,d;$ we will denote
this possibility "inside independence version". \par

 If in addition {\it all} the r.v. $ \{ \xi(i,m) \} $ are also independent, we will call
this possibility "common (total) independent version (case)".\par

{\bf Remark 1.1} In contradiction: the random vectors $ \vec{\xi}(m) $ may be particularly
coincide.\par

\vspace{3mm}

 Further, let
 $ I = I(d) = I(d,n) = \{ (i_1, i_2, \ldots, i_d )\} $  be the set
 of indices of the form

 $$
\{ (i_1, i_2, \ldots, i_d )\}: 1 \le i_1 < i_2 < i_3 \ldots < i_{d-1} < i_d \le n,
 $$

 $ J = J(d) = J(d,n) = \{ (i_1, i_2, \ldots, i_d )\} $  be the set
 of indices of the form

 $$
J(d,n) = \{ (i_1, i_2, \ldots, n )\}: 1 \le i_1 < i_2 < i_3 \ldots < i_{d-1} \le n-1 < i_d = n,
 $$
 $ b(I) = b(i_1,i_2,\ldots,i_d) $ be a $ d - $ dimensional numerical non-random sequence,

 $$
 \xi(I) \stackrel{def}{=} \prod_{m=1}^d \xi(i_m,m),
  \xi(J) \stackrel{def}{=} \prod_{m=1}^{d-1} \xi(i_m,m), \
  \sigma^2(i,m) = {\Var} ( \xi(i,m)),
 $$

$$
Q_d = Q(d,n, \{\xi(\cdot, \cdot) \} ) = Q(d,n,b) = \sum_{I \in I(d,n)} b(I) \xi(I)
$$
be a homogeneous polynomial (random polynomial) of power $  d $ on the variables
$ \{ \xi(i,m) \}, \ m \le d $ without diagonal members; on the other  words, multiple
stochastic integral on the discrete martingales measures, martingale transform,
$ n = 1,2,\ldots,\infty; $ in the case $ n=\infty \ Q(d,\infty) $ should be understood as
a limit $ Q(d,\infty) = \lim_{n \to \infty} Q(d,n) $ if there exists with probability one. \par
 Note that

 $$
 \Var Q(d,n,b) = \sum_{I \subset I(d,n)} b^2(I) \prod_{i=1}^d \sigma^2(i_m,m);
 $$
hence if

$$
\sum_{I \subset I(d,\infty)} b^2(I) \prod_{i=1}^d \sigma^2(i_m,m) < \infty
$$
then $ Q(d,\infty) $ exists. In particular, if

$$
\sup_{i,m} \Var \xi(i,m) = \sup_{i,m} \sigma^2(i,m) < \infty, \ \sum_{I \subset I(d,\infty)} b^2(I) < \infty
$$
this condition is satisfied. \par
 We will denote further

 $$
B = B(d,n) = \{b(I): \ \sum_{I \subset I(d,n)} b^2(I) = 1 \}, \ n \le \infty,
 $$
and suppose $ b \in B = B(d,n). $ Obviously, if $ b \in B, \ \sigma^2(i,m)=1, $ then
$ \Var Q_d = 1. $\par

\vspace{4mm}
{\bf  The aim of this paper is to obtain uniform tail and moment estimation  for
polynomial martingales. }\par

\vspace{4mm}

The letter $ C, C_k, C_l(\cdot) k=1,2,\ldots $ with or without subscript  will denote a
finite positive non essential constants, not necessarily  the same at each appearance. \par

\vspace{4mm}

  The applications of these estimations in the theory of probability distributions in
Banach spaces,  method Monte-Carlo and statistics are described, e.g.,  in  \cite{Billingsley1},
\cite{Buldygin1}, \cite{Frolov1}, \cite{Grigorjeva1}, \cite{Mushtary1}, \cite{Ledoux1},
\cite{Ostrovsky2}, chapter 4, \cite{Ostrovsky5},\cite{Ostrovsky6},
\cite{Vakhania1}, chapter 5 etc. \par
 The one-dimensional case $ d=1$ see, e.g. in
 \cite{Bahr1}, \cite{Braverman2}, \cite{Carothers1}, \cite{Hitczenko1},\cite{Ibragimov1},
\cite{Johnson1}, \cite{Kwapien1}, \cite{Latala1}, \cite{Litvak1}, \cite{Ostrovsky5},
 \cite{Ostrovsky6}, \cite{Rosenthal1}, \cite{Sharachmedov1}, \cite{Utev1}.
  The multidimensional case $ d \ge 2 $ is investigated in  \cite{Ostrovsky5},
see also reference therein, especially therein was considered the case light, or
equally exponential tails of distributions,  with correspondent {\it lower estimations}. \par
\vspace{3mm}
 The offered article may be considered as a continuation of the work \cite{Ostrovsky5},
where was obtained the exponential decreasing tail estimations for $ Q_d $  under
conditions that all the summands  (and multipliers) $ \xi(i,m) $ have also
exponential decreasing tails. \par
 {\it We suppose in this article as a rule that the r.v.  $ \xi(i,m) $
have heavy tails, for instance, which do not have moments of arbitrary order.}\par

\vspace{3mm}
  The paper is organized as follows. In the section 2 included some auxiliary facts about
  the random variables belonging to the so-called Grand Lebesgue Spaces (GLS); on the other
  words, moment spaces. In the section 3 and 4 we formulate and prove the main result: the moment
   bounds for  distributions of polynomial martingales.\par
   In the fifth section we generalize obtained results into the arbitrary centered polynomials
from totally independent variables.  The $ 6^{th} $  section is devoted to the tails estimations
for multiple martingale differences sums. In the next section we obtain a moment estimations
for a reverse martingales.   \par
   The last section contains some concluding remarks.\par

\vspace{3mm}
   \section{Auxiliary facts.}
\vspace{3mm}

  In order to formulate and prove the main result in this case, we recall here for reader
 convenience some facts about so-called Grand Lebesgue Spaces \cite{Fiorenza1},
 \cite{Fiorenza2}, \cite{Fiorenza3}, \cite{Fiorenza4},
  \cite{Ivaniec1}, \cite{Ivaniec2}, \cite{Liflyand1} etc.
or equally  "moment"  spaces of random variables defined on
 fixed probabilistic space $ (\Omega, \cal{A}, {\bf P}); $ more detail description
with some applications of these spaces  see
 in \cite{Buldygin1}, \cite{Kozatchenko1}, \cite{Ostrovsky2}, \cite{Ostrovsky3},
 \cite{Ostrovsky4}, \cite{Ostrovsky5}.\par

 Notice that in the articles
\cite{Fiorenza2}, \cite{Fiorenza3}, \cite{Fiorenza4},
\cite{Ivaniec1}, \cite{Ivaniec2}, \cite{Liflyand1}, \cite{Ostrovsky3}, \cite{Ostrovsky4}
was considered more general case when the measure $ {\bf P} $ may be unbounded. \par

 Let us consider the following norm (the so-called "moment norm")
on the set of r.v. defined in our probability space by the following way: the space
$ G(\nu) = G(\nu;r) $ consist, by definition, on all the  r.v. with finite norm

$$
||\xi||G(\nu) \stackrel{def}{=} \sup_{p \in (2,r)} [|\xi|_p/\nu(p)], \ |\xi|_p :=
{\bf E}^{1/p} |\xi|^p.\eqno(2.1)
$$
  Here $ r = \const > 2, \ \nu(\cdot) $ is some continuous positive on the
 {\it semi-open} interval  $ [1,r) $ function such that

     $$
     \inf_{p \in (2,r)} \nu(p) > 0, \ \nu(p) = \infty, \ p > r;
     $$
 and as usually

 $$
 |\xi|_p \stackrel{def}{=}  \left[{\bf E} |\xi|^p \right]^{1/p}
 $$

We will denote

$$
 \supp (\nu) \stackrel{def}{=} [2,r) = \{p: \nu(p) < \infty \}
$$
in the case when $ \nu(r)= \infty $ or correspondingly if $ \psi(r) < \infty $

$$
 \supp (\nu) \stackrel{def}{=} [2,r] = \{p: \nu(p) < \infty \}.
$$

 The case $ r = +\infty $ is investigated in  \cite{Kozatchenko1}, \cite{Ostrovsky5},
 \cite{Ostrovsky2}; therefore, we suppose further $  2 < r < \infty. $\par

 Let $  \xi $ be a r.v. such that there is a value $r, 2 < r < \infty,  $ for which
 $$
 \forall p > r \Rightarrow
 |\xi|_p \stackrel{def}{=} \left[ {\bf E}|\xi|^p \right]^{1/p} = \infty.
 $$

  We can and will assume furthermore that

  $$
  r = r_{\xi} =  \inf \{p: \ \left[ | \xi |_p = \infty \right]  \}=
 \sup \{p: \left[ | \xi |_p < \infty \right] \}.
  $$

  In order to  emphasis the dependance of a functions on the variable $ r, $ we will denote the
class of such a functions as $ G(\nu; r). $\par

 The {\it natural} function   $  \nu_{\xi}(p)  $ may be defined as follows:

$$
\nu_{\xi}(p) := |\xi|_p = \left[ {\bf E}|\xi|^p \right]^{1/p}.\eqno(2.2)
$$
 Obviously,

 $$
 ||\xi||G(\nu_{\xi})=1.
 $$

 The {\it natural} function for the {\it family} of a random  variables  $ \{\xi(\cdot) \}, $
or equally for a random vector $ \vec{\xi} = \{\xi(\cdot) \}
 \ \{\nu(\cdot) \} = \nu_{ \{ \xi(k) \} }(p) =  \nu_{ \vec{\xi}}(p),
\ k = 1,2,\ldots, \ n \le \infty  $ may be defined as follows:

$$
\{\nu_{ \vec{\xi}}(p) \} =
\nu_{ \{\xi\}}(p) := | \{ \vec{\xi } \}|_p = \sup_k \left[ {\bf E}|\xi(k)|^p \right]^{1/p},\eqno(2.2a)
$$
if there exists  and is finite  for the values $ p \in [1,r), \ r > 1. $ \par

 {\bf Remark 2.0.} From the Lyapunov's inequality follows that if for the r.v. $ \xi  \
|\xi|_2 < \infty,$ then $ \forall p \in [1,2] \ |\xi|_p \le |\xi|_2 = C < \infty. $\par

\vspace{4mm}

 The  complete description of a possible natural functions see in
\cite{Ostrovsky2}, \cite{Ostrovsky5}, chapter 1, section 3.\par

{\bf Example 2.1.} Suppose the r.v. $ \xi $ satisfies the condition

 $$
 T(x) = T_{\xi}(x) = x^{-r} \ \log^{\gamma}(x) \ L(\log x),  r = \const > 0,
 \ x > e, \eqno(2.3)
 $$
where as ordinary the tail function $ T(x) = T_{\xi}(x)$ for the r.v. $  \xi $ may be
defined as follows:

$$
T(x) = T_{\xi}(x) = {\bf P}(|\xi| \ge x), \ x > 0; \eqno(2.4)
$$
Here $ L = L(x) $ is positive continuous slowly varying as $ x \to \infty $ function,
$ r = \const > 1 $ and $ \gamma > -1. $ Then (see \cite{Liflyand1}, \cite{Ostrovsky2}) for
the values $ p \in [1,r), \ p \to r-0 $

$$
{\bf E}|\xi|^p \sim C(r,\gamma,L) \ (r-p)^{-\gamma-1} \ L(1/(r-p)). \eqno(2.5)
$$

 Therefore,  for $ 1 \le p < r $

 $$
\nu_{ \{\xi\}}(p) \asymp  (r-p)^{-(\gamma-1)/p} \ L^{1/p}(1/(r-p)) \asymp
(r-p)^{-(\gamma-1)/r} \ L^{1/r}(1/(r-p)).
 $$
  Note that if $ \gamma < -1 $ or more generally
 $$
 \int_e^{\infty} x^{-1} \ \log^{\gamma}(x) \ L(\log x) \ dx < \infty,
 $$
then
$$
 {\bf E} |\xi|^p < \infty \ \Leftrightarrow p \in [1,r].
$$
 Recall that we consider only the case $ p \ge 1. $ \par

 This condition is equivalent to the follows:

 $$
 \nu_{ \{\xi\} }(p) \asymp 1, \ p \in [1, r]; \ \nu_{ \{\xi\} }(p)= \infty, p > r.
 $$

  Note that an inequality $ \exists r > 0 \ \Rightarrow \psi(r-0) < \infty $ is equivalent to
the moment restriction $  |\xi|_r < \infty. $  \par

\vspace{3mm}
 Notice that the consistent estimation of regular varying distributions is investigated, e.g. in
\cite{Davydov1}, \cite{Davydov2}, \cite{Davydov3}, \cite{Paulauskas1}. The properties of slowly
and regular varying functions see in \cite{Seneta1}; see also
\cite{Bingham1}.

\vspace{3mm}
 We recall now the relations between moments for r.v. $  \xi, \ \xi \in G(\nu,r) $
 and its tail behavior. Namely, for $  p < r $

 $$
 |\xi|_p = \left[ p \int_0^{\infty} u^{p-1} \ T_{|\xi|}(u) \ du \right]^{1/p},
 $$
therefore

$$
||\xi||G(\nu;r) = \sup_{p < r}
\left\{\left[ p \int_0^{\infty} u^{p-1} \ T_{|\xi|}(u) \ du \right]^{1/p}/\nu(p) \right\}.
$$
 Conversely, if the r.v. $  \xi $ belongs to the space $ G(\nu,r), $ then

$$
T_{|\xi|}(x) \le T^{(||\xi||G(\nu;r) \cdot \nu)}(x),
$$
where by definition
$$
T^{(\nu)}(x) \stackrel{def}{=}
\inf_{p \in (1,r)} \left[ \nu^p(p)\ /x^p \right], \ x > 0.
$$

 Let us denote for simplicity  the union of all the classes $ G(\nu,r), \ r > 2 $
 with finite values $ r $ greater  as  two as  $ G(\nu):  $
 $$
G(\nu)= \cup_{2 < r < \infty}G(\nu,r).
 $$

{\bf Example 2.2.} If

$$
{\bf E}|\xi|^p  \le C_2 \ (r-p)^{-\gamma-1} \ L(1/(r-p)), \ p < r,
$$
then

$$
T_{|\xi|}(x) \le C_4(r,\gamma,L) \ x^{-r} \ [\log x]^{\gamma+1} \ L(\log x), \ x > e.
$$
 Notice that there is  a "logarithmic gap" between  upper and lower tail and moment relations.
 This gap is essential, see \cite{Liflyand1}, \cite{Ostrovsky2}. Let us consider the following \par
 {\bf Example 2.3.} Let $ \zeta $ be a discrete r.v. with distribution

$$
{\bf P} \left( \zeta = \exp(e^k)  \right) = C_5 \ \exp \left(\beta r k - r  e^k \right), \
$$
$$
 k = 1,2,\ldots, \ \beta = \const > 0,
$$
where obviously

$$
1/C_5 = \sum_{k=1}^{\infty} \exp \left(\beta r k - r  e^k \right).
$$
 We conclude after some calculations:

 $$
 |\zeta|_p \asymp (r-p)^{-\beta},
 $$
but for the sequence $ x(k)= \exp(\exp(k)) $
$$
T_{\zeta}(x(k)) \ge c_6 [\log x(k)]^{\beta} \ (x(k))^{-r}.
$$

\vspace{3mm}
{\bf Definition 2.1.} \par
\vspace{3mm}
 Let $ \nu_1 \in G(\nu;r_1), \ \nu_2 \in G(\nu;r_2), r_1,r_2 > 1  $ and let
 $ r = r(p_1,p_2):= [1/r_1 + 1/r_2]^{-1} > 1. $  We define a {\it new operation}
 $ \nu := \nu_1\otimes \nu_2 $ as follows:

$$
\nu(p) = \nu_1\otimes \nu_2(p) \stackrel{def}{=}
\inf \left\{ \nu_1 \left(\frac{p}{a}  \right) \cdot
\nu_2 \left( \frac{p}{b} \right): \ a,b >0, a + b = 1  \right\}. \eqno(2.6)
$$
 Note that the function $ \nu $ is correctly defined on the set $ 1 \le p < r. $ \par
  The following properties of operation $ \otimes $ are obvious:

  $$
  \nu_1\otimes \nu_2(p) = \nu_2\otimes \nu_1(p);
  $$

 $$
 \left(C_1  \nu_1 \right) \otimes \left( C_2  \nu_2 \right)(p) =
 C_1 C_2 \ \nu_1 \otimes \nu_2(p);
 $$

 $$
 \left(C_1 p^{\alpha} \nu_1 \right) \otimes \left( C_2 p^{\beta} \nu_2 \right)(p) \le
 C_3 \ p^{\alpha+\beta} \ \nu_1 \otimes \nu_2(p);
 $$

 $$
 \nu_1\otimes (\nu_2 + \nu_3)(p) \le \nu_1\otimes \nu_2(p) + \nu_1\otimes \nu_2(p).
 $$

  The sense of this function: as follows from the H\"older inequality,

$$
||\xi \eta||G(\nu,r) \le ||\xi||G(\nu_1,r_2) \ ||\eta||G(\nu_2, r_2). \eqno(2.7)
$$
 The {\it iteration} of this operation may be defined recursively:
 if $ \nu_m(\cdot) \in G(\nu, r_m) $ and

 $$
 r:= \left(\sum_{m=1}^d 1/r_m \right)^{-1} > 1,
 $$
 then
  $$
 \otimes_{m=1}^d \nu_m (p) =  ((\nu_1 \otimes \nu_2) \ldots \otimes \nu_d(p)), \ 1 \le p < r, \eqno(2.8)
 $$
 and \par

 {\bf Proposition 2.1.}

 $$
 ||\prod_{m=1}^d \eta(m)||G(\otimes_{m=1}^d \nu_m) \le
 \prod_{m=1}^d ||\eta(m)||G(\nu_m).
 $$

 Further, if the r.v. $ \eta(m) $ are common (total) independent  and $ s = \min_m r_m > 1,$ then

 $$
 |\prod_{m=1}^d \eta(m)|_p = \prod_{m=1}^d |\eta(m)|_p, \ 1 \le p < s;
 $$
therefore in this case

 $$
 ||\prod_{m=1}^d \eta(m)||G(\prod_{m=1}^d \nu_m,s)  =
 \prod_{m=1}^d \left[ ||\eta(m)||G(\nu_m,r_m) \right]. \eqno(2.9)
 $$
 Obviously, the product $ \prod_{m=1}^d \nu_m(p) $ is infinite iff $ p > s. $ \par
\vspace{3mm}
{\bf Example 2.4.} Let

$$
\nu_1(p) = \psi_1(r_1-p), \ \nu_2(p) = \psi_2(r_2-p),
$$
where $ r = (1/r_1+1/r_2)^{-1} > 1 $ and $ \nu_j(p) = + \infty, \ p > r_j, \ j=1,2; $
so that the point $ z=0+ $ is unique point of singularity for the functions $ \psi_m(z). $ \par
 We get after simple calculations:

$$
\nu_1 \otimes \nu_2(p) \le \psi_1 \left(\frac{r_1+r_2}{r_2} (r-p)\right) \cdot
\psi_2 \left(\frac{r_1+r_2}{r_1} (r-p)\right), \ p < r. \eqno(2.10)
$$

 More generally, let

$$
\nu_m(p) = \psi_m(r_m-p), \ m-1,2,\ldots,d,
$$

$$
  r:= (1/r_1+1/r_2+ \ldots 1/r_d)^{-1} > 1; \psi_j(p) = + \infty, \ p > r_j, \ j=1,2,\ldots,d.
$$
Denote

$$
z_j = \frac{\prod_{s\ne j}r_s }{\sum_{m=1}^d \prod_{l \ne m}r_l}.
$$
 We have using proposition (2.1)  denoting

 $$
 \psi(p) = \prod_{m=1}^d \psi_m(z_m (r-p)): \eqno(2.11)
 $$

$$
||\prod_{m=1}^d \eta(m)||G(\psi,r) \le \prod_{m=1}^d ||\eta(m)||G(\psi_m,r_m). \eqno(2.12)
$$

 If the functions $ \psi_m(\cdot)(p) $ satisfy the  $ \Delta_2 $ condition
at the points $ p = 0+: $

$$
\sup_{p \in [1,r_m)} \frac{ \psi_m(z_m (r_m-p))}{ \psi_m(r_m-p)} < \infty, \eqno(2.13)
$$
 then the function $ \psi = \psi(p) $ in (2.11) may be replaced by following expression:

 $$
 \psi(p)\asymp \overline{\psi}(p) \stackrel{def}{=} \prod_{m=1}^d \psi_m(r-p).
 $$
The condition (2.13) is satisfied, e.g., for the function of a view

$$
\psi_m(p) = C \ (r_m-p)^{-\Delta_m} \ L_m(1/(r_m-p)), \ \Delta_m = \const > 0, \eqno(2.14)
$$
$ L_m(x) $ are slowly varying as $ x \to \infty $ functions. \par
 We obtain in the considered case

 $$
 \psi(p) \asymp (r-p)^{-\sum_{k=1}^d \Delta_k} \ \prod_{k=1}^d L_k(1/(r-p)). \eqno(2.15)
 $$
More detail, if

$$
|\eta(m)|_p \le C_m \  (r_m-p)^{-\Delta_m} \ L_m(1/(r_m-p)), \ \Delta_m = \const > 0, 1 \le p < r_m,
\eqno(2.16)
$$
then

$$
\left|\prod_{m=1}^d \eta(m) \right|_p \le C^{(d)} \cdot (r-p)^{-\sum_{k=1}^d \Delta_k} \cdot
\prod_{k=1}^d L_k(1/(r-p)),\ 1 \le p < r. \eqno(2.17)
$$
 Recall that $ r = [\sum_{m=1}^d 1/r_m]^{-1} > 1.  $\par

\vspace{3mm}

{\bf  Example 2.5.} We prove here the precision of inequality (2.17) by  means of construction of
corresponding example.  Let $ \epsilon $ be a positive r.v. defined on some sufficiently rich
probabilistic space, such that

$$
T_{\epsilon}(x) = {\bf P}(\epsilon > x) = 1/x, \ x > 1.
$$
 Put
 $$
 \xi = \epsilon^{1/r_1}, \ \eta = \epsilon^{1/r_2}, \ r_1, r_2 > 1,
  r:= (1/r_1+ 1/r_2)^{-1} > 1.
 $$
We have:

$$
|\xi|^p_p \asymp \int_1^{\infty} x^{p-1-r_1} \ dx = (r_1-p)^{-1};
$$

$$
|\xi|_p \sim (r_1-p)^{-1/r_1}, \ 1 \le p < r_1,
$$
and analogously
$$
|\eta|_p \sim (r_2-p)^{-1/r_2}, \ 1 \le p < r_2.
$$

 Further, $ \xi \ \eta = \epsilon^{1/r},   $

$$
|\xi \ \eta|_p \sim (r-p)^{-1/r} = (r-p)^{-1/r_1 - 1/r_2} , \ 1 \le p < r,
$$
i.e. there is an asymptotical equality in the inequality (2.17). \par

\vspace{3mm}

{\bf Example 2.6.}  Let $ r_1 = r_2 = \infty, $

$$
|\xi|_p \le C_1 \ p^{\mu_1}, \ |\eta|_p \le C_2 \ p^{\mu_2}, \ \mu_1,\mu_2 = \const > 0,
$$
i.e. $ \nu_1(p) = C_1 \ p^{\mu_1}, \ \nu_2(p) = C_2 \ p^{\mu_2}. $  \par
 Recall, see \cite{Ostrovsky5}, \cite{Ostrovsky2}, that this means

 $$
 T_{\xi}(x) \le \exp \left(-C_3 x^{1/\mu_1} \right), \
 T_{\eta}(x) \le \exp \left(-C_4 x^{1/\mu_2} \right).
 $$
 Subexamples: in the case $ \mu_1 = 1/2 $ the r.v. $ \xi $ is {\it subgaussian:}

 $$
 T_{\xi}(x) \le C_1 \exp(-C_2 x^2), \ x \ge 0.
 $$

In the case $ \mu_1 = 1 $ the r.v. $ \xi $ is {\it subexponential:}

 $$
 T_{\xi}(x) \le C_3 \exp(-C_4 x), \ x \ge 0.
 $$
 The r.v. $ \xi $ is said to be {\it pre-gaussian,} if $ \exists \ \lambda_0,
  0 < \lambda_0 \le \infty \ \Rightarrow  $

$$
\forall \lambda, \ |\lambda| < \lambda_0 \ \Rightarrow {\bf E}\exp(\lambda \ \xi) < \infty.
$$
 Evidently, if $ \xi $ is pre-gaussian then
 $$
 \forall p \ge 0 \ \Rightarrow \nu_{\xi}(p) < \infty .
 $$

\vspace{4mm}
Further, we have:

$$
\nu_1 \otimes \nu_2(p) = C_1 \ C_2 \ \mu_1^{-\mu_1} \ \mu_2^{-\mu_2} \ (\mu_2+\mu_2)^{\mu_1+ \mu_2} \
p^{\mu_1 + \mu_2}, \eqno(2.18)
$$
or equally

$$
T_{\xi \eta}(x) \le \exp \left(-C_5 x^{1/(\mu_1 + \mu_2)}   \right). \eqno(2.18a)
$$
 The inequality (2.18) is exact, for instance, for independent variables $ \xi, \eta, $
 for which

 $$
 T_{\xi}(x) = \exp \left(-C_3 x^{1/\mu_1} \right), \
 T_{\eta}(x) = \exp \left(-C_4 x^{1/\mu_2} \right). \eqno(2.19)
 $$
In detail, it follows from (2.19) as $ p \to \infty $

$$
|\xi|_p \sim C_6 p^{\mu_1}, \ |\eta|_p \sim C_7 p^{\mu_2},
$$
and consequently

$$
|\xi \eta|_p \sim C_8 \ p^{\mu_1 + \mu_2}.
$$

 For instance, the product of two subgaussian (dependent, in general case) variables is
sub-exponential.\par

\vspace{3mm}

{\bf Example 2.7.} Let now $ r = \const > 1, $

$$
|\xi|_p \le C_1 (r-p)^{-\gamma} \ L_1(1/(r-p)) =: \nu_1(p), \ 1 \le p < r,
$$

$$
|\eta|_p \le C_2 p^{\mu} L_2(p) =:\nu_2(p), \ p \ge 1,
$$
where $ L_1(z), L_2(z) $  are slowly varying as $ z \to \infty $  positive continuous functions.\par
 We obtain substituting in the formula (2.6) the values $ a=1/(1+C(r-p) ), \ b = 1-a: $

 $$
 \nu_1 \otimes \nu_2(p) =  |\xi \ \eta|_p \le C_3 (r-p)^{-\gamma - m} \ L_1(1/(r-p)) \
 L_2(1/(r-p)), \ 1\le p < r. \eqno(2.20)
 $$
  We show now the exactness of an estimation (2.20). Let $ \Omega =  \{\omega \} =[0,1] $ with
Lebesgue measure $ {\bf P} $  and let

$$
\xi = \omega^{-1/r} \ |\log \omega|^{\kappa} \ L(|\log \omega|), \ r = \const > 1,
$$
 where $ L = L(x) $ is as before is positive continuous slowly varying as $ x \to \infty $ function,
$$
  \eta = |\log \omega|^{\mu}, \ \mu = \const > 0.
$$
 We find by the direct calculation, see \cite{Ostrovsky3}:

 $$
 |\xi|_p \asymp (r-p)^{- \kappa-1/r} \ L(1/(r-p)), p \in [1,r);
 \ |\eta|_p \asymp p^{\mu}, \ p \ge 1;
 $$

$$
|\xi \ \eta|_p \asymp (r-p)^{-\mu - \kappa -1/r} \ L(1/(r-p)), p \in [1,r).
$$

\vspace{3mm}

{\bf Example 2.8.}  Define the following $ \nu $ function  $ \psi_{(r)}(p):$

$$
\psi_{(r)}(p) = 1, \ p \le r; \ \psi_{(r)}(p)= +\infty, p > r, \ r = \const > 1.
$$
 It is easy to verify using the Lyapunov's inequality that the $ G\psi_{(r)}(\cdot) $  space
coincides with the classical Lebesgue-Riesz space $ L(r): $

$$
||\xi|| G\psi_{(r)}(\cdot) = |\xi|_r. \eqno(2.21)
$$
 Further, let $ \xi(i), \ i=1,2,\ldots,d $ be some r.v. from the spaces $  G\psi_{(r(i))} = L(r(i)), $
where as before $ r(i) > 1, \ r:= 1/(\sum 1/r(i)) \ge 1. $  We infer:

$$
||\prod_i \xi(i) ||G\psi_{(r)} \le \prod_i ||\xi(i)||G\psi_{(r(i))}, \eqno(2.22a)
$$
or equally

$$
|\prod_i \xi(i)|L(r) \le \prod_i |\xi(i)|L(r(i)).\eqno(2.22b)
$$
 For instance, $ |\xi \eta|_1 \le |\xi|_2 \cdot |\eta|_2; \ |\xi \eta \zeta|_1 \le
|\xi|_3 \cdot |\eta|_3 \cdot \zeta|_3. $ \par
 Evidently, the estimates (2.22a) or (2.22b) are essentially non-improvable. \par

\vspace{3mm}

\section{Main result: martingale version}

\vspace{3mm}

 We recall before formulating the main result  some useful for us
moment inequalities  for the  sums of centered  martingale differences
and independent r.v., \cite{Ostrovsky4}. Namely, let $ \{ \theta(i) \} $ be a
sequence of centered  martingale differences relative any filtration; then

$$
\sup_n \sup_{b \in B} \left|\sum_{i=1}^n b(i) \theta(i) \right|_p \le K_M(p) \
\sup_i |\theta(i)|_p, \eqno(3.0)
$$
where for the {\it optimal value} of the constant $ K_M = K_M(p) $ is true the inequality

$$
 K_M(p) \le   p \ \sqrt{2}, \ p \ge 2.
$$
 Note that the upper bound in (3.0)
 $$
  K_I(p) \le  0.87 p/\log p, \ p \ge 2
 $$
is true for the independent centered r.v. $ \{\theta(i)) \}, $ see also \cite{Ostrovsky4}. \par
\vspace{3mm}
 We denote the natural function for the sequences $ \vec{\xi}(m) =
 \{ \xi(i,m) \}, i=1,2,\ldots,n $ as $ \nu_{m}= \nu_{m}(p): $

$$
\nu_{m}(p) = \sup_{k=1,2,\ldots,n} |\xi(k,m)|_p, \ m=1,2,\ldots,d, \eqno(3.1)
$$
and suppose that for all the values $ 1 \le m \le d $ the function $ \nu_{m}(\cdot) $
belongs to the set $ G(\nu).$  More detail, we assume

$$
\nu_{m}(\cdot) \in G(\nu, r_m), \ r_m > 2.
$$

 Define the sequence of a functions $ \zeta_m = \zeta_m(p), \ m=1,2,\ldots,d  $
by the following recursion:

$$
\zeta_1(p) = K_M(p) \ \nu_1(p), \eqno(3.2)
$$
(initial condition) and

$$
\zeta_{m+1}(p) =  K_M(p) \cdot \left[ \left\{ \zeta_m \otimes \nu_{m+1} \right\}(p) \right],
\ m=1,2,\ldots, d-1. \eqno(3.3)
$$

\vspace{3mm}

{\bf Theorem 3.1 (martingale version).} \par
 Let us denote as before
 $$
 r = \left(\sum_{m=1}^d \frac{1}{r_m} \right)^{-1} \eqno(3.4)
 $$
and suppose $ r > 1. $ Proposition:
$$
\sup_{b \in B}||Q(d,n,b)|| G (\zeta_{d}) \le
\prod_{m=1}^d  ||\vec{\xi}(m)|| G(\nu^{(m)}, r_m). \eqno(3.5)
$$

\vspace{3mm}

{\bf Proof} is at the same as in \cite{Ostrovsky5}; it used the estimate (3.0)
and induction over $ d. $ Namely, the one-dimensional case $ d=1 $ coincides directly
with inequality (3.0). \par

 Let now $ d \ge 2. $ The sequence $ (Q_d(n), F(n)) $ is again martingale with
 correspondent martingale-differences

 $$
 \beta(n) := Q_d(n)-Q_d(n-1) = \xi(n,d) \sum_{I \in J(d,n) } b(I)\xi(I)=
 $$

$$
\xi(n,d) \sum_{1 \le i_1 < i_2 \ldots < i_{d-1} \le n-1}b(i_1,i_2,\ldots,i_{d-1}, n)
\prod_{m=1}^{d-1} \xi(i_m,m) =
$$

$$
\left[\xi(n,d) \cdot \sqrt{\sum_{I \in J(d,n)} b^2(I) }  \right] \times
$$
$$
\left[ \frac{\sum_{I \in J(d,n)} b(I)\xi(I)}{ \sqrt{\sum_{I\in J(d,n)} b^2(I)} } \right]
\stackrel{def}{=} \{ \eta(n,d) \} \cdot \{\tau(n,d)\}, \eqno(3.6)
$$
where

$$
\tau(n,d) := \left[ \frac{\sum_{I \in J(d,n)} b(I)\xi(I)}{ \sqrt{\sum_{I\in J(d,n)} b^2(I)} } \right],
$$

$$
\eta(n,d) := \left[\xi(n,d) \cdot \sqrt{\sum_{I \in J(d,n)} b^2(I) }  \right].\eqno(3.7)
$$

It follows from the induction statement

$$
|\tau(n,d)|_p \le  \zeta_{d-1} (p).\eqno(3.8)
$$

 It remains to use  the induction statement and the definition of the operation
$ \otimes: $\par

$$
|Q_d(n)|_p \le \nu_{d} \otimes  \zeta_{d-1} (p), \eqno(3.9)
$$

\vspace{3mm}

{\bf Example 3.1.} Let

$$
\nu_{m}(p) \le (r_m-p)^{-\Delta_m} \ L_m(1/(r_m-p)), \ 2 <r_m < \infty, \
\Delta_m = \const > 0,
$$

$$
r:= (\sum_{m=1}^d 1/r_m)^{-1} > 1,  \ 1 \le p < r.
$$
$ L_m(x) $ be positive continuous slowly varying as $ x \to \infty $ functions. \par
 As long as  a multiplier $ p \ \sqrt{2} $ is bounded in considered case:
$ p \ \sqrt{2} \le r \ \sqrt{2}, $  we conclude by virtue of theorem 3.1:

$$
\sup_{b \in B}||Q(d,n,b)|| G (\zeta_{d}) \le
\prod_{m=1}^d  ||\vec{\xi}(m)|| G(\nu_{m}, r_m),
$$
where

$$
\zeta_{d}(p) = 2^{d/2} \ r^{d}  \ (r-p)^{-\sum_{m=1}^d \Delta_m} \
\prod_{m=1}^d  L_m(z_m/(r-p)),
$$
or equally

$$
\sup_{b \in B} |Q_d|_p \le 2^{d/2} \ r^{d}  \ (r-p)^{-\sum_{m=1}^d \Delta_m} \
\prod_{m=1}^d  L_m(z_m/(r-p)), \ 1 \le p < r. \eqno(3.10)
$$

\vspace{3mm}

{\bf Example 3.2.} As long as the case when $ \forall m \ r_m = \infty $ was
investigated in  \cite{Ostrovsky5},  we consider here  the {\it mixed} possibility,
when

$$
\exists k \in [2,d-1], \ r_1,r_2,\ldots, r_k \in (2, \infty); \ r_{k+1} = r_{k+2} =
\ldots = r_d = \infty.
$$
 In detail, suppose

 $$
 \nu_{m}(p) = (r_m - p)^{-\Delta_m} L_m(1/(r_m-p)), \ m=1,2,\ldots,k, \ 2 \le p < r_m;
 $$

$$
\nu_{m}(p) = p^{\mu_k} \ L_m(p), \ m = k+1,\ldots,d, \ p \in [1, \infty),
$$
where as before  all the functions $ L_m(x) $ be positive continuous slowly varying as
$ x \to \infty $ functions. \par
 Put

 $$
 \overline{r} = (\sum_{m=1}^k 1/r_m)^{-1}>1, \ \overline{L}(p) =
 \prod_{m=1}^k L_m(z_m/(\overline{r}-p)) \ \prod_{m=k+1}^d L_m(1/(\overline{r}-p), \ 1 \le p < \overline{r}.
 $$

 From theorem 3.1 follows:

$$
\zeta_{d}(p) = 2^{k/2} \ \overline{r}^{k}  \ ( \overline{r}-p)^{-\sum_{m=1}^k \Delta_m \ -\sum_{m=k+1}^d \mu_m
-(d-k)} \ \overline{L}(p). \eqno(3.11)
$$

\vspace{3mm}

{\bf Example 3.3.}  Let us illustrate the {\it exactness} of the assertions of examples (3.1)
and (3.2). It is sufficient thereto to consider the case $ n=1, d=2 $  and refer to the examples
2.4, 2.5 and 2.6.\par

\vspace{3mm}
{\bf Example 3.4.}  Let  the source centered martingale differences
$ \xi(i,m), \ i=1,2,\ldots,n $ be  r.v. from
the fixed finite ball of the  spaces $  G\psi_{(r(i))} = L(r(i)): $

$$
\sup_i |\xi(i,m)|_{r(m)} =: D(m) < \infty,
$$
where as before $ r(i) > 1, \ r:= 1/(\sum_m 1/r(m)) \ge 1. $  We deduce
by virtue of the example 2.8 and the boundedness of the value $ p \sqrt{2}: $

$$
\sup_{b \in B} |Q_d|_r \le C_3(d; \{r(m)\}) \ \prod_{m=1,2, \ldots,d}
\sup_{i=1,2,\ldots,n} |\xi(i,m)|_{r(m)}. \eqno(3.12)
$$

  Obviously, the estimate (3.12) is unprovable; it is sufficient to show this to consider
the case $ d=2,3, \ldots; \ n=1 $ and refer to the example 2.8. \par

\vspace{3mm}

\section{Main result: independent versions}

\vspace{3mm}

 Analogously may be considered  all the "independent cases". \par

\vspace{3mm}

{\bf A. "Common independent version". } \par

\vspace{3mm}

 Recall  that in this case all the
random variables $ \{ \xi(i,m) \}, i=1,2,\ldots,n; \ m=1,2,\ldots,d $
are mean zero and independent.\par

\vspace{3mm}

 {\bf Theorem 4.1.}\par

 Define the sequence of a functions $ \zeta_m = \zeta_m(p), \ m=1,2,\ldots,d  $
by the following way:

$$
\zeta_1(p) = K_I(p) \ \nu_1(p), \eqno(4.1)
$$
(initial condition) and  by recursion

$$
\zeta_{m+1}(p) =  K_M(p) \ \zeta_m(p) \times \nu_{m+1}(p), \ m=1,2,\ldots, d-1; \eqno(4.2)
$$
with the explicit solution

$$
\zeta_d(p) = K_I(p) \ K_M^{d-1}(p) \ \prod_{m=1}^d \nu_m(p). \eqno(4.3)
$$

Proposition:
$$
\sup_{b \in B}||Q(d,n,b)|| G (\zeta_{d}) \le
\prod_{m=1}^d  ||\vec{\xi}(m)|| G(\nu_{m}, r_m). \eqno(4.4)
$$

\vspace{3mm}

{\bf Example 4.1.} Let the all the r.v. $ \xi(i,m) $ are mean zero and common (total)
 independent, $ \xi(i,m) \in G(\nu_m, r_m), $  such that

$$
\sup_{i=1,2,\ldots,n} ||\xi(i,m)||G(\nu_m, r_m) < \infty,
$$
and
$ \max_m r_m < \infty, \ s = \min_m r_m > 1,$ then

 $$
  \zeta_d(p) \le C  \prod_{m=1}^d \sup_{i=1,2,\ldots,n} |\xi(i,m)|_p, \ 1 \le p < s;
 $$
therefore in this case

 $$
 \sup_{b \in B} ||Q_d||G(\prod_{m=1}^d \nu_m,s)  \le C \cdot
 \prod_{m=1}^d  \sup_{i=1,2,\ldots,n} ||\xi(i,m)||G(\nu_m, r_m).
  $$
 Evidently, the product $ \prod_{m=1}^d \nu_m(p) $ is infinite iff $ p > s. $ \par

\vspace{3mm}

{\bf B. "Inside independent version". } \par

\vspace{3mm}

{\bf Theorem 4.2. }.\par

\vspace{3mm}

 Define the sequence of a functions $ \zeta_m = \zeta_m(p), \ m=1,2,\ldots,d  $
by the following way:

$$
\zeta_1(p) = K_I(p) \ \nu_1(p), \eqno(4.5)
$$
(initial condition) and

$$
\zeta_{m+1}(p) =  K_M(p) \ \left[ \zeta_m(p) \otimes \nu_{m+1}(p) \right], \ m=1,2,\ldots, d-1; \eqno(4.6)
$$
(recursion). Proposition:

$$
\sup_{b \in B}||Q(d,n,b)|| G (\zeta_{d}) \le
\prod_{m=1}^d  ||\vec{\xi}(m)|| G(\nu_{m}, r_m). \eqno(4.7)
$$

\vspace{3mm}

{\bf Example 4.2} is a  slight modification of Example 3.1. Indeed, let under condition of the
theorem 4.2

$$
\nu_{m}(p) \le (r_m-p)^{-\Delta_m} \ L_m(1/(r_m-p)), \ 2 <r_m < \infty, \
\Delta_m = \const > 0,
$$

$$
r:= (\sum_{m=1}^d 1/r_m)^{-1} > 1,  \ 1 \le p < r.
$$
$ L_m(x) $ be positive continuous slowly varying as $ x \to \infty $ functions. \par
 As long as  a multiplier $ p \ \sqrt{2} $ is bounded in considered case:
$ p \ \sqrt{2} \le r \ \sqrt{2}, $  we conclude by virtue of theorem 3.1:

$$
\sup_{b \in B}||Q(d,n,b)|| G (\zeta_{d}) \le
 C( d, \{ r_m \}) \ \prod_{m=1}^d  ||\vec{\xi}(m)|| G(\nu_{m}, r_m),
$$
where

$$
\zeta_{d}(p) = \ (r-p)^{-\sum_{m=1}^d \Delta_m} \ \prod_{m=1}^d  L_m(z_m/(r-p)),
$$
or equally

$$
\sup_{b \in B} |Q_d|_p \le C( d, \{ r_m \})  \ (r-p)^{-\sum_{m=1}^d \Delta_m} \
\prod_{m=1}^d  L_m(z_m/(r-p)), \ 1 \le p < r.
$$

\bigskip

{\bf C. "Vector independent version".} \par

\vspace{3mm}

{\bf Theorem 4.3. }.\par

\vspace{3mm}

 Define the sequence of a functions $ \zeta_m = \zeta_m(p), \ m=1,2,\ldots,d  $
by the following way:

$$
\zeta_1(p) = K_M(p) \ \nu_1(p), \eqno(4.8)
$$
(initial condition) and

$$
\zeta_{m+1}(p) =  K_M(p) \ \zeta_m(p) \times \nu_{m+1}(p),
 \ m=1,2,\ldots, d-1; \eqno(4.9)
$$
(recursion); with the explicit solution:

$$
\zeta_d(p) = K_M^d(p) \cdot \prod_{m=1}^d \nu_m(p),
$$

Proposition:

$$
\sup_{b \in B}||Q(d,n,b)|| G (\zeta_{d}) \le
\prod_{m=1}^d  ||\vec{\xi}(m)|| G(\nu_{m}, r_m). \eqno(4.10)
$$
\vspace{3mm}

\section{Arbitrary centered polynomial from independent variables}

\vspace{3mm}

{\it We generalize here the results of the last section on the
arbitrary centered polynomial from common independent variables. }\par
\vspace{3mm}
 In detail: let all the random variables $ \{\xi(i,m) \}, \ 1=1,2,\ldots,n; \ m=1,2,\ldots,d $
be centered and common independent. Arbitrary centered polynomial $ R_d =R_d(n) = R_d(n,\{\xi(i,m) \} $
of degree $ d $ may be written in the form of

$$
R_d = \sum_{1 \le i_1 < i_2 \ldots < i_d \le n }
b( i_1,i_1,\ldots,i_1; i_2,i_2,\ldots,i_2; \ldots,i_d,i_d, \ldots,i_d ) \times
$$

$$
\prod_{l=1}^{d} \left[ \xi^{k(l)}(i_l,l) - m(k(l),i_l,l) \right], \eqno(5.1)
$$
where $ m(k,i,l) = {\bf E}\xi^k(i,l), \ k(l) \in \{0,1,\ldots,d  \}, $  so that
$ \sum_l k(l)=d, $

$$
k(l) = \card \{i_l  \ {\bf in} \ b( i_1,i_1,\ldots,i_1; i_2,i_2,\ldots,i_2; \ldots,i_d,i_d, \ldots,i_d ) \},
$$
and if some $ k(l)=0 $ then by definition

$$
\xi^{k(l)}(i_l,l) - m(k(l),i_l,l) = 1.
$$

 We again suppose $ b \in B;$  the multi-sequence $ b(I) $ may be not symmetric. \par
 The  convergence of the series for $ R_d $ is investigate in \cite{Ostrovsky5}.\par
 For instance, if $ d=2, $ then $ R_d $ has a view:

 $$
 R_d = \sum \sum_{1 \le i < j \le n} b(i,j) \xi(i) \xi(j) + \sum_{i=1}^n b(i,i)
 (\xi^2(i)- {\bf E}\xi^2(i)) +
 $$

 $$
  \sum \sum_{1 \le i < j \le n} c(i,j) \eta(i) \eta(j) + \sum_{i=1}^n c(i,i)
 (\xi^2(i)- {\bf E}\xi^2(i)) +
 $$
 $$
 \sum \sum_{1 \le i < j} a_1(i,j) \xi(i) \eta(j) + \sum \sum_{1 \le i < j} a_2(i,j) \eta(i) \xi(j).
 $$

{\bf Theorem 5.1.} Suppose that

$$
T_{\xi(i,m)}(x) \le x^{-r(m)} \ \log^{\gamma(m)}(x) \ L_m(\log x), \ x \ge e, \eqno(5.2)
$$
where $  r(m), \gamma(m) = \const $ such that $ r(m) < \infty, \ L_m(z)  $ are continuous slowly
varying as $ z \to \infty $ positive functions,

$$
\underline{r} \stackrel{def}{=} \min_m r(m) >d.
$$
We define:

$$
M(r) = \{m: \ r(m)= \underline{r} \}, \ \overline{\gamma} = \max_{m \in M(r)} \gamma(m),
$$

$$
M_{\gamma} = \{m: \ m \in M(r), \gamma(m) = \overline{\gamma} \},
$$

$$
\overline{L}(z) = \max_{m \in M_{\gamma}} L_m(z), z \ge 1.
$$
 Proposition: for the values  $ p, 1 \le p < \underline{r}/d $ and in the case when
$  \gamma > -1 $

 $$
\sup_{b \in B} {\bf E} |R_d|^p \le C( \{r(m)\}, \{\gamma(m)\}, \{L_m(\cdot)\})\cdot
 (\underline{r}/d-p)^{-\gamma-1} \ \overline{L}(1/(\underline{r}/d-p)). \eqno(5.3a)
 $$
\vspace{3mm}

 When $ \gamma < - 1 $ or more generally

 $$
\max_{m_0 \in M_{\gamma} } \sup_{j=1,2, \ldots,n} {\bf E} |\xi^d(j,m_0)|^p < \infty,
 $$
then

$$
\sup_{b \in B} {\bf E} |R_d|^p \le C_1 ( \{r(m)\}, \{\gamma(m) \} ) < \infty. \eqno(5.3b)
$$
\vspace{3mm}
{\bf Proof.} Let $ m_0 $ be arbitrary number from the set  $ M_{\gamma}; $ it is easy to verify
that the main member in the expression for $ R_d $ has a view:

$$
R_M = \sum_{j=1}^n b(j) ( \xi^d(j,m_0) - {\bf E}\xi^d(j,m_0)), \eqno(5.4)
$$
where $ \sum_j b^2(j)=1. $ \par

 We have from the condition (5.2)

 $$
 T_{ \xi^d(j,m_0) - {\bf E}\xi^d(j,m_0)}(x) \le
C \  x^{- \underline{r}/d} \ \log^{\overline{\gamma}}(x) \ \overline{L}(\log x), \ x \ge e, \eqno(5.5)
 $$
therefore

$$
{\bf E}|\xi^d(j,m_0) - {\bf E}\xi^d(j,m_0)|^p \le C_2 \
(\underline{r}/d-p)^{-\gamma-1} \ \overline{L}(1/((\underline{r}/d-p))).\eqno(5.6)
$$
 It remains to use the Rosenthal's inequality for the r.v. $ R_M. $\par
Note that in connection with common independence

$$
|\xi(1,i_1) \xi(2,i_2) \ldots \xi(d,i_d)|_p = |\xi(1,i_1)|_p |\xi(2,i_2)|_p \ldots |\xi(d,i_d)|_p  \le
$$

$$
C_3 \ (r(1)-p)^{-\gamma(1)-1} \ L_1((r(1)-p)) \ (r(2)-p)^{-\gamma(2)-1} \ L_2((r(2)-p))) \ldots \times
$$
$$
(r(d)-p)^{-\gamma(d)-1} \ L_d((r(d)-p))), \ p < \min_m r(m), \eqno(5.7)
$$
which is significantly smallest as the right-hand side of inequality (5.6). \par
\vspace{3mm}

{\bf Remark 5.1.} The assertion of theorem 5.1. may be rewritten as follows. Let us denote

$$
\rho_d(p) = (\underline{r}/d-p)^{-\gamma-1} \ \overline{L}(1/(\underline{r}/d-p)), \ 1 \le p <
\underline{r}/d;
$$

$$
\theta_m(p) = (r(m)-p)^{-\gamma(m)-1} \ L_m(1/(r(m)-p)), \ 1 \le p < r(m),
$$

then

$$
\sup_{b \in B} ||R_d||G(\rho_d) \le C( \{r(m)\}, \{\gamma(m)\}, \{L_m(\cdot)\})
\prod_{m=1}^d \ \left[\sup_{i=1}^n ||\xi(i,m)||G(\theta_m) \right]. \eqno(5.8)
$$

 Note that our results improve and generalize ones obtained in \cite{Schudy1}. \par

\vspace{3mm}

\section{Tail estimations }

\vspace{3mm}

 In many practical cases is convenient to operate by tails estimations for sums of
random variables instead moments ones, for example in the method Monte-Carlo and statistics,
\cite{Ostrovsky6},  see also \cite{Braverman2}. \par

 Let $ \nu = \nu(p) $ be any function from the set $ G(\nu), $ for instance, from the set $ G(\nu,r);$
we denote

$$
T^{(\nu)}(x) = \exp \left[ - \left(p \log \nu(p) \right)^*(\log x) \right], \ x > e, \eqno(6.1)
$$
where

$$
f^*(y) = \sup_{z, z > 0, f(z) < \infty} (yz - f(z)) \eqno(6.2)
$$
is ordinary Young-Fenchel, or Legendre transform of the function $ f = f(z). $  \par

 Further, let $ \zeta_d = \zeta_d(p) $ (or $ \rho_d = \rho_d(p)) $ be one of the functions
$ \zeta_d(p), \ \rho_d(p) $  introduced before. \par

{\bf Theorem 6.1.}  We have correspondingly for the martingale, independent and polynomial
cases  and for the values $ x > e: $

$$
\sup_{b \in B} T_{Q_d}(x) \le T^{(\zeta_d)}(x/C_1), \eqno(6.3a),
$$

$$
\sup_{b \in B} T_{V_d}(x) \le T^{(\zeta_d)}(x/C_2), \eqno(6.3b),
$$

$$
\sup_{b \in B} T_{R_d}(x) \le T^{(\rho_d)}(x/C_3), \eqno(6.3c).
$$
{\bf Proof} follows immediately from the obtained moment estimations and from the following
implication \cite{Liflyand1}, \cite{Ostrovsky3}: if for the r.v. $ \eta $ there holds
$ ||\eta||G(\nu) \le 1 $ or equally

$$
|\eta|_p \le \nu(p), \ p: \nu(p) < \infty,
$$
then

$$
T_{\eta}(x) \le T^{(\nu)}(x), \ x > e. \eqno(6.4)
$$

\vspace{3mm}

{\bf Example 6.1.} If for some $  r > 1 $

$$
\zeta_d^p(p)  \le C \ (r-p)^{-\gamma-1} \ L(1/(r-p)), \ 1 \le p < r,
$$
where as before $ L = L(z) $ be positive continuous slowly varying as $ z \to \infty $ function,
then

$$
\sup_{b \in B} T_{Q_d}(x) \le C_4(r,\gamma,L) \ x^{-r} \ [\log x]^{\gamma+1} \ L(\log x), \ x > e.
\eqno(6.5)
$$

\vspace{3mm}

{\bf Example 6.2.}  Let  the source centered martingale differences
$ \xi(i,m), \ i=1,2,\ldots,n $ be  r.v. from
the fixed finite ball of the  spaces $  G\psi_{(r(i))} = L(r(i)): $

$$
\sup_i |\xi(i,m)|_{r(m)} =: D(m) < \infty,
$$
where as before $ r(i) > 1, \ r:= 1/(\sum_m 1/r(m)) \ge 1. $  We deduce
using the Tchebychev's inequality and inequality (3.12)

$$
\sup_{b \in B} {\bf P}(|Q_d|>x) \le C_3(d; \{r(m)\}) \ \prod_{m=1,2, \ldots,d}
\sup_{i=1,2,\ldots,n} |\xi(i,m)|_{r(m)} \cdot x^{-r}, \ x > 0. \eqno(6.6)
$$

\vspace{3mm}

\section{Reverse martingales }

\vspace{3mm}

 We investigate {\it in this section} the case when all the sequences $ \vec{\xi}(m)=
 \{ \xi(i,m) \}, \ m=1,2,\ldots,d; \ i=n,n+1,\ldots,N, $ where $ 1 \le N \le \infty $ are
mean zero reverse martingales differences {\it  relative some reverse filtration } $  F(i). $ \par
 This means that

$$
 F(i+1) \subset F(i) \subset F,
$$
and $  \forall k=i+1, i+2, \ldots  \ \Rightarrow $
$$
{\bf E}\xi(i,m)/F(k)=0; \ {\bf E} \xi(i,m)/F(i)= \xi(i,m) \ (\mod {\bf P}), \eqno(7.1)
$$
"reverse martingale version".\par
 Analogously to the first section may be defined "vector independent version", "inside
independent version" and "common independent version"  for reverse martingales. \par
 An example: let $ \xi(i,m) $ be a sequence of independent centered r.v. Set

 $$
 S_m(n) = \sum_{i=n}^N \xi(i,m), \ F(i) = \sigma(\xi(i,m), \xi(i+1,m), \ldots), \eqno(7.2)
 $$
where we suppose in the case $ N=\infty $ that the series in (7.2) converge with probability
one and in the $ L_2(\Omega) $ sense.\par
 For all the values $ m $ the pair $ (S_m(n), F(n)) $ is a reverse martingale. The correspondent
reverse martingales differences are

$$
S_m(n) - S_m(n+1) = \xi(n,m).
$$

  More facts about martingales and reverse martingales see in the classical book \cite{Hall1}.
For instance, for the reverse martingale $ (S(n), F(n)) $ there exists a limit
$ \lim_{n \to \infty} S(n) $ with probability one.  Further, let $ N < \infty. $ \par
\vspace{3mm}
{\bf Proposition 7.1.}\par
\vspace{3mm}
The sequence $ (Z_i, F(i)), \ i=1,2,\ldots,N $ is reverse martingale iff the sequence

$$
(Z_{N-i+1}, F_{N-i+1})  \eqno(7.3)
$$
is ordinary martingale, see  \cite{Hall1}, chapter 1, section 1. \par
\vspace{3mm}

Let $ I=I_n^{(N)} = I_n^{(N)}(d) = \{ (i_1, i_2, \ldots, i_d )\} $  be the set
 of indices of the form

 $$
\{ (i_1, i_2, \ldots, i_d ) \}: n \le i_1 < i_{2} <  \ldots < i_{d-1} < i_d \le N,
 $$

 $ J = J_n^{(N)}(d) = \{ (i_1, i_2, \ldots, i_d )\} $  be the set
 of indices of the form

 $$
J_n^{(N)}(d) = \{ (i_1, i_2, \ldots, N )\}: n \le i_1 < i_2 < i_3 \ldots < i_{d-1} \le N-1 < i_d = N,
 $$
 $ b(I) = b(i_1,i_2,\ldots,i_d) $ be a $ d - $ dimensional numerical non-random sequence,

 $$
 \xi(I) \stackrel{def}{=} \prod_{m=1}^d \xi(i_m,m),
  \xi(J) \stackrel{def}{=} \prod_{m=1}^{d-1} \xi(i_m,m), \
  \sigma^2(i,m) = {\Var} (\xi(i,m)),
 $$

$$
V_d = V_n^{(N)}(d,n, \{\xi(\cdot, \cdot) \} ) = \sum_{I \in I_n^{(N)}(d,n)} b(I) \xi(I) \eqno(7.4)
$$
be a homogeneous polynomial (random polynomial) of power $  d $ on the variables
$ \{ \xi(i,m) \}, \ m \le d $ without diagonal members.

 Let us denote
$$
\nu_m^{(n)}(p) = \sup_{i: n \le i \le N} |\xi(i,m)|_p,\eqno(7.5)
$$
where in the case $ N < \infty $ the symbol $"\sup" $ must be replaced by $"\min";$

$$
\zeta_d(p) = \zeta_d^{(n)}(p) = \sup_{b \in B} | V_n^{(N)}(d,n, \{\xi(\cdot, \cdot) \} )|_p =
 \sup_{b \in B}|V_d|_p. \eqno(7.6)
$$

 {\bf Remark 7.1.} It will be presumed that for some $ b(m), 2 < b(m) \le \infty $
the value $ \nu_m^{(n)}(p) $ is finite for $ 2 \le p < b(m).  $\par

 {\bf Remark 7.2.}It is naturally {\it to wait} that in the case $ N = \infty \ \forall m= 1,2,\ldots,d \
 \Rightarrow  \lim_{n \to \infty} \sum_{i=n}^N \xi(i,m) = 0 $ and hence

 $$
 \lim_{n \to \infty}\nu_m^{(n)}(p) =0, \ 2 \le p < b(m),
$$
or equally that there exist a sequences $ \epsilon_n^{(m)} $ tending to null as $ n \to \infty $
such that

$$
\nu_m^{(n)}(p) \le  \epsilon_n^{(m)}\cdot  \nu_m^{(0)}(p), \eqno(7.7)
$$
where the sequences $  \nu_m^{(0)}(p) $ are bounded for all the values $p, \ 2 \le p < b(m).  $\par
\vspace{3mm}

{\bf Theorem 7.1. (Reverse martingale case.)} \par
We define the sequence  $  \zeta_d(p) $  by the following "reverse" recurrent relation

 $$
\zeta_m(p) = K_M(p) \cdot \left\{ \left[ \nu_m^{(n)}\otimes \zeta_{m+1}\right](p) \right\}, \
m=d-1,d-2, \ldots,1 \eqno(7.8)
 $$
with the following endpoint condition:

$$
\zeta_d(p) = K_M(p)\cdot \nu_d^{(n)}(p). \eqno(7.9)
$$

Let us denote as before
 $$
 r = \left(\sum_{m=1}^d \frac{1}{r_m} \right)^{-1}
 $$
and suppose $ r > 1. $ Proposition:
$$
\sup_{b \in B}|| V_n^{(N)}(d,n, \{\xi(\cdot, \cdot) \} ) || G (\zeta_{1}) \le
\prod_{m=1}^d  ||\vec{\xi}(m)|| G(\nu_m^{(n)}). \eqno(7.10)
$$

\vspace{3mm}

{\bf Theorem 7.2 ( "Total independent version").}\par

 Define the sequence of a functions $ \zeta_m = \zeta_m(p), \ m=1,2,\ldots,d  $
by the following way:

$$
\zeta_d(p) = K_I(p) \ \nu_d^{(n)}(p), \eqno(7.11)
$$
(endpoint condition) and  by reverse recursion

$$
\zeta_{m}(p) =  K_M(p) \ \zeta_{m+1}(p) \times \nu_{m}^{(n)}(p), m=d-1,d-2, \ldots,1 \eqno(7.12)
$$
with the explicit solution

$$
\zeta_1(p) = K_I(p) \ K_M^{d-1}(p) \ \prod_{m=1}^d \nu_m(p). \eqno(7.13)
$$

Proposition:
$$
\sup_{b \in B} || V_n^{(N)}(d,n, \{\xi(\cdot, \cdot) \} ) || G (\zeta_{1}) \le
\prod_{m=1}^d  ||\vec{\xi}(m)|| G(\nu^{(m)}, r_m). \eqno(7.14)
$$

\vspace{3mm}

{\bf Theorem 7.3 ("Inside independent version").}\par

\vspace{3mm}

 Define the sequence of a functions $ \zeta_m = \zeta_m(p), \ m=1,2,\ldots,d  $
by the following way:

$$
\zeta_d(p) = K_I(p) \ \nu^{(n)}_d(p), \eqno(7.15)
$$
(endpoint condition) and

$$
\zeta_{m}(p) =  K_M(p) \cdot \left\{ \left[ \zeta_{m+1} \otimes \nu_{m+1}^{(n)}(p) \right](p) \right\}, \
m=d-1, d-2, \ldots,1 \eqno(7.16)
$$
( reverse recursion). Proposition:

$$
\sup_{b \in B} || V_n^{(N)}(d,n, \{\xi(\cdot, \cdot) \} ) || G (\zeta_{1}) \le
\prod_{m=1}^d  ||\vec{\xi}(m)|| G(\nu_{m}^{(n)}, r_m). \eqno(7.17)
$$

\vspace{3mm}

{\bf Theorem 7.4 ("Vector independent version"). } \par

\vspace{3mm}

 Define the sequence of a functions $ \zeta_m = \zeta_m(p), \ m=1,2,\ldots,d  $
by the following way:

$$
\zeta_d(p) = K_M(p) \ \nu_d^{(n)}(p), \eqno(7.18)
$$
(endpoint condition) and

$$
\zeta_{m}(p) =  K_M(p) \ \zeta_{m+1}(p) \times \nu^{(n)}_{m}(p),
 \ m=d-1,d-2,\ldots, 1; \eqno(7.19)
$$
(reverse recursion); with the explicit solution:

$$
\zeta_1(p) = K_M^d(p) \cdot \prod_{m=1}^d \nu_m^{(n)}(p).
$$

 Proposition:

$$
\sup_{b \in B} || V_n^{(N)}(d,n, \{\xi(\cdot, \cdot) \} ) || G (\zeta_{1}) \le
\prod_{m=1}^d  ||\vec{\xi}(m)|| G(\nu_{m}^{(n)}, r_m). \eqno(7.20)
$$

\vspace{3mm}

{\bf Proofs} are at the same as in the sections 3 and 4 and may omitted. \par

\vspace{3mm}

\section{Concluding remarks}

\vspace{3mm}

{\bf 1. A "good $ \lambda $  inequality".} \par

\vspace{3mm}

 Let $ X $ and $ Y $ be two  non-negative r.v. such that there exist
a constants $ \beta > 1, \ \delta > 0, \ \epsilon > 0  $  for which

$$
\forall \lambda > 0 \ {\bf P}(X > \beta \lambda, \ Y \le \delta \lambda)
\le {\bf P}(X > \lambda).
$$
 The following inequality (implication) is called {\it A "good $ \lambda $  inequality":}

$$
{\bf E} X^p \le \beta^p \delta^{-p} (1-\epsilon \beta^p)^{-1} {\bf E } Y^p, \
0 < p < |\log_{\beta} \epsilon|.
$$
 This assertion is proved, e.g. in \cite{Hall1}, chapter 2, lemma 2.4 and is used
widely in the theory of martingales.\par
 Let us denote $ r= |\log_{\beta} \epsilon| $ and suppose $ r > 1, $

$$
Y \in G(\psi,r) \ \Leftrightarrow |Y|_p \le \psi(p), \ p < r.
$$
 For the function

 $$
 \psi_r(p) := \psi(p) \ (r-p)^{-1/r}
 $$
we deduce:

$$
||X||G(\psi_r,r) \le C \ ||Y||G(\psi,r).
$$

\vspace{3mm}

{\bf 2. Maximal inequality for polynomial martingales.} \par

\vspace{3mm}

 As long as all the sequences $ n \to Q(d,n,b) $ and $ n \to V^{(N)}(d,N-n,b) $ are
martingales relative the source filtration $  F(i) $ (or $ F(N-i)), $  we can apply the
famous Doob's inequality under before formulated conditions:

$$
\sup_{b \in B} \left| \sup_{n \le N} Q(d,n,b) \right|_p  \le C \zeta_d(p) \cdot \frac{p}{p-1},
$$

$$
\sup_{b \in B} \left| \sup_{n \le N} V^{(N)}(d,n,b) \right|_p  \le C \zeta_d(p) \cdot \frac{p}{p-1}.
$$

\vspace{3mm}

{\bf 3. A new norm of a random vectors.}  \par

\vspace{3mm}
 Let $ X = X(\omega)= \vec{X} = (X_1,X_2,\ldots, X_D) $
be a $ D- $ dimensional random vector and $ b \in S^{D-1} $  be
arbitrary non-random normed vector  with at the same dimension as the vector $ X: $
$$
b \in S^{D-1} \Leftrightarrow b \in B \Leftrightarrow \sum_{i=1}^D b^2(i) = 1.
$$
 We define the $ L(p), \ p \ge 1 $ and correspondingly $ G(\psi) $ norm of the random
vector (vector $ L_p(\Omega) $ norm) $ X = \vec{X} = (X_1,X_2,\ldots, X_D) $  as follows:

$$
|\vec{X}|_p \stackrel{def}{=} \sup_{b \in B} |(X,b)|_p,
$$

$$
||\vec{X}||G(\psi) = ||\vec{X}||G(\psi;a,b)
\stackrel{def}{=} \sup_{p \in (a,b)} [|\vec{X}|_p/\psi(p)].
$$
 Note that in the one-dimensional case $ D=1 $ this definition coincides with ordinary. \par
If for example  all the components $ X_i $ are centered, independent and identical distributed,
then

$$
|X_1|_p \le |\vec{X}|_p \le K_I(p) \ |X_1|_p.
$$

 Our main results, for instance, theorem 3.1, may be reformulated in the terms of vector
 $ L_p $  and vector $ G(\psi) $ norms  as follows:

$$
||\xi(I)|| G (\zeta_{d}) \le
\prod_{m=1}^d  ||\vec{\xi}(m)|| G(\nu^{(m)}, r_m) = 1.
$$
 Here

 $$
 D = \frac{n(n-1)(n-2)\ldots (n-d+1)}{d!}, \ n > d+1.
 $$

\vspace{3mm}

{\bf 4. Rosenthal's approach.}  \par

\vspace{3mm}

 H.P.Rosenthal in  \cite{Rosenthal1} proved that each sequence of centered independent r.v.
generates in $ L_p(\Omega), \ 2 \le p < \infty  $ the subspace isomorphic to $ l_2. $ See
also  \cite{Braverman2}. \par
 Our main results imply that this is true for the multiple sequences of  mean zero polynomial
martingales $ \xi(I): $

$$
|b(\cdot)|^2 = \sum_I b^2(I).
$$

\vspace{3mm}

{\bf 5. Normed multiple sums estimates.} \par

\vspace{3mm}
 It is no hard to obtain the moment and tail estimates for  {\it normed} multiple sums

$$
\hat{Q_d} = \frac{Q_d}{ \sqrt{\Var (\{Q_d\})} }
$$
in the terms of {\it normed} summands

$$
\hat{\xi}(i,m) = \frac{\xi(i,m)}{ \sqrt{\Var\{\xi(i,m) \}} } =
\frac{\xi(i,m)}{ \sigma(i,m) }.
$$
 They are at the same as before; when we write instead the functions $ \nu_m(p) $
the functions

$$
\hat{\nu}_m(p) = \sup_{i=1,2,\ldots,n}|\hat{\xi}(i,m)|_p =
\sup_{i=1,2,\ldots,n} \left[ \frac{\xi(i,m)}{\sigma(i,m)} \right]_p
$$
and correspondingly $ b \in \hat{B}, $ where

$$
\hat{B} = \{b = \vec{b}: \sum_{I \subset I(d,n) } b^2(I) \prod_{m=1}^d \sigma^2(i_m,m) = 1 \}.
$$

\end{document}